\documentclass[11pt, reqno]{amsart}
\allowdisplaybreaks[3]

\usepackage[indentafter]{titlesec}
\usepackage{enumerate}
\usepackage{mathrsfs}
\usepackage{amsopn,amscd,amsmath,amssymb,amsthm,amsfonts,epsfig,graphics}
\usepackage{pslatex,avant,color}
\usepackage{cite}

\newcommand{\ran}{{\rm ran}}

\newcommand{\dom}{{\rm dom}}
\newcommand{\mul}{{\rm mul}}

\theoremstyle{plain}
\newtheorem{thm}{Theorem}[section]
\newtheorem{lem}{Lemma}[section]
\newtheorem{cor}{Corollary}[section]
\newtheorem{propo}{Proposition}[section]
\newtheorem{defi}{Definition}[section]
\newtheorem{exm}{Example}[section]

\theoremstyle{remark}

\titleformat{name=\section}{}{\thetitle.}{0.8em}{\centering\scshape}
\titleformat{name=\subsection}{}{\thetitle.}{0.5em}{\bfseries}[]

\newcounter{str}

\numberwithin{equation}{section}
\begin{document}
 \vskip .8cm

\title[Characterization of Weyl functions]{Characterization of Weyl functions in the class of operator-valued generalized Nevanlinna functions}
\author{Muhamed Borogovac}

\begin{abstract}
We provide the necessary and sufficient conditions for a generalized Nevanlinna function $Q$ ($Q\in N_{\kappa }\left( \mathcal{H} \right)$) to be a Weyl function (also known as a Weyl-Titchmarch function).   

We also investigate an important subclass of $N_{\kappa }(\mathcal{H})$, the functions that have a boundedly invertible derivative at infinity $Q'\left( \infty \right):=\lim \limits_{z \to \infty}{zQ(z)}$. These functions are regular and have the operator representation $Q\left( z \right)=\tilde{\Gamma}^{+}\left( A-z \right)^{-1}\tilde{\Gamma},z\in \rho \left( A \right)$,  where $A$ is a bounded self-adjoint operator in a Pontryagin space $\mathcal{K}$. We prove that every such strict function $Q$ is a Weyl function associated with the symmetric operator $S:=A_{\vert (I-P)\mathcal{K}}$, where $P$ is the orthogonal projection, $P:=\tilde{\Gamma} \left( \tilde{\Gamma}^{+} \tilde{\Gamma} \right)^{-1} \tilde{\Gamma}^{+} $. 

Additionally, we provide the relation matrices of the adjoint relation $S^{+}$ of $S$, and of $\hat{A}$, where $\hat{A}$ is the representing relation of $\hat{Q}:=-Q^{-1}$. We illustrate our results through examples, wherein we begin with a given function $Q\in N_{\kappa }\left( \mathcal{H} \right)$ and proceed to determine the closed symmetric linear relation $S$ and the boundary triple $\Pi$ so that $Q$ becomes the Weyl function associated with $\Pi$. 
\end{abstract}


\maketitle
\thispagestyle{empty}

Dedicated to Prof. Mirjana Vukovic for her jubilee.

\section{Introduction}\label{s2}
1.1. We denote the sets of positive integers, real numbers, and complex numbers by $\mathbb{N}$, $\mathbb{R}$, and $\mathbb{C}$, respectively. Let $\left( \mathcal{K},\, \left[ .,. \right] \right)$ represent a Krein space. That is a complex vector space equipped with a scalar product $\left[ .,. \right]$, which is a Hermitian sesquilinear form. It admits the following decomposition of $\mathcal{K}$:
\[
\mathcal{K}=\mathcal{K}_{+}[+]\mathcal{K}_{-},
\]
where $\left( \mathcal{K}_{+},[.,.] \right)$ and $\left( \mathcal{K}_{-},-[.,.] \right)$ are 
Hilbert spaces that are mutually orthogonal with respect to the form $[.,.]$. Elements $x,y \in \mathcal{K} $ are \textit{orthogonal} if $\left[ x,y\right]=0$, denoted by $x \left[ \perp \right]y$. Every Krein space $\left( \mathcal{K},[.,.] \right)$ is \textit{associated} with a Hilbert space $\left( \mathcal{K},(.,.) \right)$, defined as a direct and orthogonal sum of the Hilbert spaces $\left( \mathcal{K}_{+},[.,.] \right)$ and $\left( \mathcal{K}_{-},-[.,.] \right)$. The topology in the Krein space $\mathcal{K}$ is induced by the associated Hilbert space $\left( \mathcal{K},(.,.) \right)$. The \textit{orthogonal companion} $A^{\left[ \perp \right]}$  of the set $A$ is defined by $A^{\left[ \perp \right]}:=\left\lbrace y\in \mathcal{K}:x\left[ \perp \right]y,\forall x \in A\right\rbrace $, and the \textit{isotropic} part $M$ of $A$ is defined by  $M:= A \cap A^{\left[ \perp \right]}$. For properties of Krein spaces, one can refer to e.g., \cite[Chapter V]{Bog}.

If the scalar product $\left[ .,. \right]$ has $\kappa \in \mathbb{N}$ negative squares, then we call it a \textit{Pontryagin space of negative index} $\kappa $. If $\kappa =0$, then it is a Hilbert space. More information about Pontryagin space can be found, for example, in \cite {IKL}. 

The following definitions of a linear relation and basic concepts related to it can be found in \cite {A,S,DM2}. In the following, $X$, $Y$, and $W$ represent Krein spaces which include Pontryagin and Hilbert spaces. 

A \textit{linear relation} $T: X\rightarrow Y$ is a linear manifold $T \subseteq X \times Y $.

If $ X=Y $, then $T$ is said to be a \textit{linear relation in} $X$. A linear relation $T$ is closed if it is a (closed) subspace with respect to the product topology of $X \times Y $. As usual, for a linear relation or operator $ T: X \rightarrow Y $, or $ T \subseteq X \times Y $, the symbols $ \dom T $, $ \ran\, T $, and $ \ker T $ represent the domain, range and kernel, respectively. Additionally, we will use the following concepts and notation for two linear relations, $T$ and $S$ from $X$ into $Y$, and a linear relation $U$ from $Y$ into $W$:

\[
\mul\, T:=\left\{ g\in Y : \left\{ 0,g \right\}\in T \right\},
\]
\[
T\left( f \right):=\left\{ g\in Y, : \left\{ f,g \right\}\in T \right\}, (f\in D\left( T \right)),
\]
\[
T^{-1}:=\left\{ \left\{ g,f \right\}\in Y \times X : \left\{ f,g \right\}\in T \right\},
\]
\[
zT:=\left\{ \left\{ f,zg \right\}\in X \times Y : \left\{ f,g \right\}\in T \right\}, (z\in \mathbb{C}),
\]
\[
S+T:=\left\{ \left\{ f,g+k \right\} : \left\{ f,g \right\}\in S,\left\{ f,k \right\}\in T \right\},
\]
\[
S\hat{+}T:=\left\{ \left\{ f+h,g+k \right\} : \, \left\{ f,g \right\}\in S,\left\{ h,k \right\}\in T \right\},
\]
\[
S\dot{+}T:=\left\{ \{f+h,g+k\} : \left\{ f,g \right\}\in S,\left\{ h,k \right\}\in T, S\cap T=\lbrace 0\rbrace \right\},
\]
\[
UT:=\left\{ \left\{ f,k \right\}\in X \times W \, : \, \left\{ f,g \right\}\in T,\left\{ g,k \right\}\in U\, for\, some\, g\in Y \right\},
\]
\[
T^{\ast}:=\left\{ \left\{ k,h \right\}\in Y \times X : \left[ f,h \right]=\left[ g,k \right]\, for\, all\, \left\{ f,g \right\}\in T \right\},
\]
\[
T_{\infty }:=\left\{ \left\{ 0,g \right\}\in T \right\}.
\]
If $T(0)=\{ 0\}$, we say that T is \textit{single-valued} linear relation, i.e.\textit{ operator}. The sets of closed linear relations, closed operators, and bounded operators in X are denoted by $\tilde{C}(X)$, $C(X)$, $B(X)$, respectively.  

Let $A$ be a linear relation in a Krein space $\mathcal{K}$. When $X=Y=\mathcal{K}$ we use the notation $A^{+}$ rather than $A^{\ast}$. We say that $A$ is \textit{symmetric} (\textit{selfadjoint}) if it satisfies $A\subseteq A^{+}$ ($A=A^{+})$.  

Every point $\alpha \in \mathbb{C}$ for which $\left\{ f,\alpha f \right\}\in A$, with some $f\ne 0$, is called a \textit{finite eigenvalue}, denoted by $\alpha \in \sigma_{p}(A)$. The corresponding vectors are \textit{eigenvectors belonging to the eigenvalue} $\alpha $. If for some $z \in \mathbb{C}$ the operator $(A-z)^{-1}$ is bounded, not necessarily densely defined in $\mathcal{K}$, then $z$ is a \textit{point of regular type of} $A$, symbolically, $z\in \hat{\rho }\left( A \right)$. If for $z\in \mathbb{C}$ the relation $\left( A-z \right)^{-1}$ is a bounded operator and $\overline{\ran \left(A-z \right)} = \mathcal{K}$, then $z$ is a \textit{regular point of} $A$, symbolically $z\in \rho \left( A \right)$.  

In a Pontryagin space $\mathcal{K}$, an isometric operator $U$ is called \textit{unitary} if $\dom\, U= \ran\, U=\mathcal{K}$, see \cite[Definition 5.4]{IKL}. 

According to the definition \cite[Definition 1.3.7]{BHS}, linear relations $T:\mathcal{K}\rightarrow\mathcal{K} $ and $T':\mathcal{K}'\rightarrow\mathcal{K}'$ are unitarily equivalent if there exists a unitary operator $U:\mathcal{K}\rightarrow\mathcal{K}'$ such that $T'=\lbrace \lbrace U(x),U(x')\rbrace:  \lbrace x,x'\rbrace\in T\rbrace$.

Let $\mathcal{L}(\mathcal{H})$ denote the Banach space of bounded operators in a Hilbert space $\mathcal{H}$. Recall that an operator valued function $Q:\mathcal{D}\left( Q \right) \subset \mathbb{C} \to \mathcal{L}(\mathcal{H})$ belongs to the \textit{generalized Nevanlinna class} $N_{\kappa }\left( \mathcal{H} \right)$ if it is meromorphic on $\mathbb{C} \backslash \mathbb{R}$, such that ${Q\left( z \right)}^{\ast }=Q\left( \bar{z} \right)$, for all points $z$ of holomorphy of $Q$, and the kernel $N_{Q}\left( z,w \right):=\frac{Q\left( z \right)-{Q\left( w \right)}^{\ast }}{z-\bar{w}}$ has $\kappa $ negative squares. A generalized Nevanlinna function $Q\in N_{\kappa }\left( \mathcal{H} \right)$ is called \textit{regular} if the operator ${Q\left( w \right)}$ is boundedly invertible at least for one point $w\in \mathcal{D}(Q)$, see \cite{Lu1}. 

We will need the following, Krein-Langer representation of generalized Nevanlinna functions.

\begin{thm}\label{theorem12}. A function $Q:\mathcal{D}(Q)\subset \mathbb{C} \to \mathcal{L}(\mathcal{H})$ is a generalized Nevanlinna function of some index $\kappa $ if and only if it has a representation of the form
\begin{equation}
\label{eq12}
Q\left( z \right)={Q(w)}^{\ast }+(z-\bar{w})\Gamma_{w}^{+}\left( I+\left( z-w \right)\left( A-z \right)^{-1} \right)\Gamma_{w}, z\in \mathcal{D}\left( Q \right),
\end{equation}
where, A is a self-adjoint linear relation in some Pontryagin space $(\mathcal{K}, [.,.])$ of index $\tilde{\kappa }\ge \kappa ; \Gamma_{w}:\mathcal{H}\to \mathcal{K}$ is a bounded operator; $w\in 
\rho \left( A \right)\cap \mathbb{C}^{+}$ is a fixed point of reference. This representation can be chosen to be minimal, that is
\begin{equation}
\label{eq14}
\mathcal{K}=c.l.s.\left\{ \Gamma_{z}h:z\in \rho \left( A \right),h\in \mathcal{H} \right\}
\end{equation}
where
\begin{equation}
\label{eq16}
\Gamma_{z}:=\left( I+\left( z-w \right)\left( A-z \right)^{-1} \right)\Gamma_{w}.
\end{equation}
If realization (\ref{eq12}) is minimal, then $\tilde{\kappa }= \kappa $. In that case $ \mathcal{D}(Q)=\rho(A)$ and the triple $(\mathcal{K},\, A,\, \Gamma _{w})$ is uniquely determined (up to unitary equivalence). 
\end{thm}

The linear relation $A$ in (\ref{eq12}) is called a \textit{representing relation (operator)} of $Q$. Such operator representations were developed by M. G. Krein and H. Langer, see e.g. \cite {KL1, KL2} and later converted to representations in terms of linear relations, see e.g. \cite{DLS, HSW}. 

Functions $Q \in N_{\kappa }(\mathcal{H})$ which fulfill the condition
\begin{equation}
\label{eq18}
\bigcap\limits_{z\in D\left( Q \right)} {\ker \frac{Q\left( z \right)-Q\left( \bar{w} \right)}{z-\bar{w}}}=\lbrace 0\rbrace
\end{equation}
for one, and hence for all, $w \in \mathcal{D}(Q)$, are called \textit{strict}, see e.g. \cite[p. 619]{BeLu}.

In what follows, $S$ denotes a closed symmetric relation or operator, not necessarily densely defined in a separable Pontryagin space $(\mathcal{K}\left[ .,. \right])$, and $S^{+}$ denotes an adjoint linear relation of $S$ in $(\mathcal{K}\left[ .,. \right])$. For definitions and notation of concepts related to an ordinary boundary triple $\Pi$ for the linear relation $S^{+}$, see e.g. \cite{BHS, D1, D2}. We copy some of those definitions here with adjusted notation. For example, the operator denoted by $\Gamma_{2}$ in \cite{D1} is denoted by $\Gamma_{0}$ in \cite{BHS,D2} and here, while $\Gamma_{1}$ denotes the same operator in all papers. Elements of $S^{+}$ are denoted by $\hat{f}, \thinspace \hat{g},\, \mathellipsis $, where e.g. $\hat{f}:= \left({\begin{array}{*{20}c}
f\\
f'\\
\end{array}}\right)=
\left\{ f,f'\right\}$. 
Let 
\[
\mathcal{R}_{z }:=\mathcal{R}_{z }(S^{+})=\ker \left( S^{+}-z \right),\, z \in \hat{\rho}(S),
\]
be the \textit{defect subspace} of $S$. Then
\begin{equation}
\label{eq122}
\hat{\mathcal{R}}_{z }:=\left\{ \left( {\begin{array}{*{20}c}
f_{z }\\
z f_{z }\\
\end{array} } \right):f_{z }\in \mathcal{R}_{z } \right\},\, \, \mathcal{R}:=\left( 
\dom\, S \right)^{\left[ \bot\right] } ,\, \, \hat{\mathcal{R}}:=\left\{ \left( 
{\begin{array}{*{20}c}
0\\
f\\
\end{array} } \right):f\in \mathcal{R} \right\}.
\end{equation}
\begin{defi}\label{definition12}. \cite[Definition 2.1]{D1} A triple $\Pi =(\mathcal{H}, \Gamma_{0}, \Gamma_{1})$, where $\mathcal{H}$ is a Hilbert space and $\Gamma_{0}, \Gamma_{1}$ are bounded operators from $S^{+}$ to $\mathcal{H}$, is called an ordinary boundary triple for the relation $S^{+}$ if the abstract Green's identity 
\begin{equation}
\label{eq124}
\left[ f',g \right]-\left[ f,g' \right]=\left( \Gamma_{1}\hat{f},\Gamma_{0}\hat{g} \right)_{\mathcal{H}}-\left( \Gamma_{0}\hat{f}, \Gamma_{1}\hat{g} \right)_{\mathcal{H}}, \forall \hat{f}, \hat{g}\in S^{+}, 
\end{equation}
holds, and the mapping $\Gamma :\hat{f}\to \left( {\begin{array}{*{20}c}
\Gamma_{0}\hat{f}\\
\Gamma_{1}\hat{f}\\
\end{array} } \right)$ from $S^{+}$ to $\mathcal{H}\times \mathcal{H}$ is surjective. 
\\
The operator $\Gamma $ is called the boundary or reduction operator.
\end{defi}


An extension $\tilde{S}$ of $S$ is called proper, if $S \subsetneq \tilde{S} \subseteq S^{+}$. The set of proper extensions of $S$ is denoted by $Ext\, S$. Two proper extensions $S_{0}, \thinspace S_{1}\in Ext\, S$ are called \textit{disjoint} if $S_{0}\cap S_{1} =\, S$, and \textit{transversal} if, additionally, 
$S_{0}\hat{+}S_{1}=\, S^{+}$ . 

Each ordinary boundary triple is naturally associated with two self-adjoint extensions 
of $S$, defined by $S_{i}:=\ker \Gamma_{i}, i = 0, 1$, i.e., we have 
$S_{i}=S_{i}^{+}, i = 0, 1$, see \cite [p. 4425]{D1}. 

Under above notation, the function 
\[
\emptyset \neq \rho (S_{0}) \ni z\mapsto \gamma_{z}=\left\lbrace \left\lbrace \Gamma_{0} \hat{f}_{z}, f_{z} \right\rbrace : \hat{f}_{z} \in  \hat{\mathcal{R}}_{z}(S^{+}) \right\rbrace 
\]
is called the $\gamma$-\textit{field} associated with the boundary triple $\Pi =(\mathcal{H}, \Gamma 
_{0}, \Gamma_{1})$, and the function 
\begin{equation}
\label{eq126}
\emptyset \neq \rho (S_{0}) \ni z\mapsto M(z)=\left\lbrace \left\lbrace \Gamma_{0} \hat{f}_{z}, \Gamma_{1} \hat{f}_{z} \right\rbrace : \hat{f}_{z} \in  \hat{\mathcal{R}}_{z}(S^{+}) \right\rbrace 
\end{equation}
is called the \textit{Weyl function} associated with the boundary triple $\Pi =(\mathcal{H}, \Gamma_{0}, \Gamma_{1})$, see e.g. \cite{DM1,BHS,D1}. Let us mention that functions $\gamma_{z}: \mathcal{H}\rightarrow \mathcal{R}_{z}$ are bijections and satisfy the formula (\ref{eq16}). 

1.2. \textbf{The following is a summary of the results presented in this paper.} Basic concepts of the Weyl function and $\gamma$-field of the symmetric operator $S$ in the Hilbert space setting were introduced in the classical papers, see \cite{DM1,DM2}. For later developments in the field of boundary relations and Weyl functions, we refer the reader to \cite{D1,DHMS1,BDHS,Be}. 

In this paper, we prove a characterization of the Weyl functions in the class of operator valued regular generalized Nevanlinna functions. Therefore, we use operator (relation) representations in the Pontryagin space $(\mathcal{K},[.,.])$ setting of the regular generalized Nevanlinna function $Q \in N_{\kappa }\left( \mathcal{H} \right)$. We denote by $A$ the representing self-adjoint relation of $Q$ and by $\hat{A}$ the representing self-adjoint relation of $\hat{Q}=-Q^{-1}$.

In Section \ref{s4}, in Proposition \ref{proposition22} and Example \ref{example22}, we show how to derive the strict part of a generalized Nevanlinna function. It is well known that a strict function need not to be invertible, see e.g. \cite{DHS}. In Example \ref{example22}, we see that a regular function $Q$ need not to be strict.

In Theorem \ref{theorem22}, one of the main results of the paper, we give a characterization of the Weyl functions in terms of regular and strict generalized Nevanlinna functions. In Theorem \ref{theorem22} (b), we prove the more difficult converse part. It is a generalization of the converse part of \cite[Theorem 1]{DM1} in several levels. Namely, in the converse part of \cite[Theorem 1]{DM1}, authors start with a Krein $\mathcal{Q}$-function of a given  symmetric operator $S$ in a Hilbert space. This means they assume the existence of the symmetric operator $S$, and then they prove the existence of the corresponding boundary triple that has the Weyl function equal to the given $\mathcal{Q}$-function. 

We solve a more general problem. We only assume that a regular and strict generalized Nevanlinna function is given, i.e. we do not assume the existence of a symmetric operator or relation $S$. We first have to prove  the existence of the symmetric linear relation $S$ in a Pontryagin space to be in a position to find the corresponding triple. In order to prove the existence of the symmetric relation $S$, we use much later results from \cite {Lu1}. 

Similar issues were studied for the definitizable matrix function, see \cite[Theorem 3.5]{Be}.

Section \ref{s6} can be viewed as an application of \cite{B} in the area of boundary triples and Weyl functions. In this section, we deal with an important subclass of regular functions $Q\in N_{\kappa }\left( \mathcal{H} \right)$, the functions that have a boundedly invertible derivative $Q'\left( \infty \right):=\lim \limits_{z \to \infty}{zQ(z)}$. We are again focused on finding a symmetric operator $S$ and a boundary triple $\Pi$ for a given function $Q$. We start with such a function $Q$ with the representing bounded operator $A$, and in Theorem \ref{theorem32} we prove that there exists a symmetric operator $S$ such that $Q$ is the Weyl function corresponding to $S$ and $A$. Hence, here we also give a solution of the converse problem. Moreover, we give matrix representations of $A$, $\hat{A}$, $S$, and $S^{+}$. Theorem \ref{theorem32} also gives us useful new relationships between linear relations $A$, $\hat{A}$, $S$, $S^{+}$ and $\hat{\mathcal {R}}$ associated with a given function $Q \in N_{\kappa}(\mathcal{H})$.

In Corollary \ref{corollary310}, we prove that $\hat A$, $A$ and $S^{+}$ are $\mathcal{R}$-regular extensions of $S$ if the corresponding function $Q$ is strict and $Q'\left( \infty \right)$ is boundedly invertible.

In Section \ref{s8}, we make use of the abstract results of sections \ref{s4} and \ref{s6}. In examples \ref{example46} and \ref{example410}, the functions have a boundedly invertible derivative $Q'(\infty)$, i.e. they satisfy the assumptions of Theorem \ref{theorem32}. Therefore, we apply Theorem \ref{theorem32} to find the closed symmetric relation $S$ and the corresponding ordinary boundary triple $\Pi$ in each of the examples so that $Q$ is the Weyl function associated with $\Pi$. In Example \ref{example410}, we use Theorem \ref{theorem32} also to find relation matrices $\hat{\mathcal {R}}$, $\hat{A}$, $S$ and $S^{+}$ for the given function $Q \in N_{\kappa}(\mathcal{H})$ represented by $A$.

In Example \ref{example48} we prove that the strict part $\tilde{Q}$ of the function $Q$ used in Example \ref{example22} is indeed a Weyl function corresponding to some symmetric relation $S$ and the corresponding boundary triple $\Pi$.

\section{Characterization of Weyl functions in the set of regular generalized Nevanlinna functions $N_{\kappa}(\mathcal{H})$}\label{s4}
2.1 We will need the following lemma and proposition.

\begin{lem}\label{lemma21} \cite[Lemma 4.2]{B1} Let $Q \in N_{\kappa}\left( \mathcal{H} \right)$ be a minimally represented function by a triplet $\left( \mathcal{K},A,\Gamma_{w}\right)$ in representation (\ref{eq12}). 
\begin{enumerate}[(i)]
\item If $z \in \mathcal{D}(Q)$, then 
\[
\ker{\Gamma_{z}}=\ker{\Gamma_{w}}=:\ker {\Gamma}; \forall w \in \mathcal{D}(Q),
\]
\item 
\[
\ker{ \Gamma}=\left\lbrace  h \in \mathcal{H}: \frac{Q\left( z \right)-Q\left( \bar{w} 
\right)}{z-\bar{w}}h = 0, \thinspace \forall z, \forall w \in \mathcal{D}(Q)\right\rbrace .
\]
\end{enumerate}
\end{lem}

According to Lemma \ref{lemma21} we can introduce the Hilbert space $\tilde{\mathcal{H}}:=\left( \ker \Gamma\right)^{\bot} $ and operators $\tilde{\gamma}_{w} := (\Gamma_{w})_{|\tilde{\mathcal{H}}}$.

\begin{propo}\label{proposition22} Let $Q \in N_{\kappa}\left( \mathcal{H} \right)$ be a function minimally represented by (\ref{eq12}) with operators $\Gamma_{z}:\mathcal{H}\rightarrow \mathcal{K}$ defined by (\ref{eq16}) that satisfy (\ref{eq14}). Then the following hold: 
\begin{enumerate}[(i)]
\item Operators $\Gamma_{z}, z\in \mathcal{D}(Q)$ are one-to-one if and only if the function $Q(z):\mathcal{H}\rightarrow \mathcal{H}$ is strict.
\item For every function $Q \in N_{\kappa}\left( \mathcal{H} \right)$ minimaly represented by (\ref{eq12}) with the triple $(\mathcal{K}, A, \Gamma_{w})$, there exists a unique, up to multiplicative constant, strict function $\tilde{Q} \in N_{\kappa}\left( \tilde{\mathcal{H}} \right)$ defined by (\ref{eq12}) with the triple $(\mathcal{K}, A, \tilde{\gamma}_{w})$. Functions $Q$ and $\tilde{Q}$ have the same number of positive squares as well. 
\end{enumerate}
\end{propo}

\begin{proof} (i) This is an obvious consequence of the previous lemma.

(ii) Since, for every $w \in \mathcal{D}(Q)=\mathcal{D}\left( \tilde{Q} \right)$, the operator $\tilde{\gamma}_{w}: \tilde{\mathcal{H}}\rightarrow \ran\, \Gamma_{w} $ coincides with $\Gamma_{w}$ everywhere except on $\ker\, \Gamma_{w}$, the Pontryagin space defined by (\ref{eq14}) with $\tilde{\gamma}_{w}$ instead $\Gamma_{w}$ coincides with $\mathcal{K}$. Because $\tilde{\gamma}_{w}, \forall w \in \mathcal{D}(\tilde{Q}) $, are injections 
\[
\bigcap\limits_{z,w \in \mathcal{D}\left( \tilde{Q} \right)} {\ker \frac{\tilde{Q}\left( z \right)-\tilde{Q}\left( \bar{w} \right)}{z-\bar{w}}}=\emptyset.
\] 
holds, i.e., $\tilde{Q}$ is a strict function. The representing relation $A$ remains the same because functions $\tilde{\gamma}_w, \forall w \in \mathcal{D}(Q)$ do not change anything in $\mathcal{K}$. 

For elements $h, k \in \mathcal{H}=\mathcal{\tilde{H}}(+)\ker \Gamma$ we have the corresponding unique orthogonal decomposition $h=\tilde{h}(+)h_{0} \wedge k=\tilde{k}(+)k_{0}$. Therefore, 
\[
\left[\frac{\tilde{Q}\left( z \right)-\tilde{Q}\left( \bar{w} \right)}{z-\bar{w}}\tilde{h}, \tilde{k} \right]=
\left[ \tilde{\gamma}_{z}\tilde{h},\tilde{\gamma}_{w}\tilde{k} \right] = \left[ \Gamma_{z}h,\Gamma_{w}k \right].
\]
This means that the numbers of both negative and positive squares of $Q$ and of $\tilde{Q}$ are the same. 
\end{proof}

The function $\tilde{Q} \in N_{\kappa}( \tilde{\mathcal{H}} )$, introduced in Proposition \ref{proposition22}, will be referred to as the \textit{strict part} of $Q$. Additionally, we will call the Hilbert space $\tilde{\mathcal{H}}$ the \textit{domain} of the strict part $\tilde{Q}$.

\begin{exm}\label{example22} Consider the following regular matrix function
\[
Q\left( z \right)=\left( {\begin{array}{*{20}c}
\frac{z}{2}-1 & \frac{z}{2}\\
\frac{z}{2} & \frac{z}{2}+1\\
\end{array} } \right).
\]
Then for vector $h=\left( {\begin{array}{*{20}c}
1\\
-1\\
\end{array} } \right)$,
\[
N(z,w)h=\frac{Q\left( z \right)-Q\left( \bar{w} \right)}{z-\bar{w}}h = 0, \forall w,z \in \mathcal{D}(Q).
\]
Therefore, this is an example of a regular function that is not strict. 

Our task is to find the strict part $\tilde{Q}$ of $Q$. 
\end{exm}
Let us switch from the basis $e_{1}= \left( {\begin{array}{*{20}c}
1\\
0\\
\end{array} } \right) $,
$e_{2}= \left( {\begin{array}{*{20}c}
0\\
1\\
\end{array} } \right) $ to the new ortho-normal basis $f_{1}= \frac{1}{\sqrt{2} }\left( {\begin{array}{*{20}c}
1\\
-1\\
\end{array} } \right) $,
$f_{2}= \frac{1}{\sqrt{2} } \left( {\begin{array}{*{20}c}
1\\
1\\
\end{array} } \right) $. With respect to the new basis, we have
\[
Q\left( z \right)=\left( {\begin{array}{*{20}c}
0 & -1\\
-1 & z\\
\end{array} } \right) \wedge f_{1}=\left( {\begin{array}{*{20}c}
1\\
0\\
\end{array} } \right) \wedge  f_{2}=\left( {\begin{array}{*{20}c}
0\\
1\\
\end{array} } \right) \wedge h=\sqrt{2}f_{1}.
\]
According to Proposition \ref{proposition22}, we can introduce the domain of $\tilde{Q}$ by $\tilde{\mathcal{H}}=l.s.\lbrace f_{2} \rbrace$. Then, if we denote by $P_{|\tilde{\mathcal{H}}}$ the orthogonal projection onto $\tilde{\mathcal{H}}$ we get the strict part of $Q$ 
\[
\tilde{Q}(z)=P_{|\tilde{\mathcal{H}}}Q(z)_{|\tilde{\mathcal{H}}}=z, z\in \mathcal{D}(Q).
\]
Recall that the strict part preserves the numbers of positive and negative squares. $\hfill$ $\square$

Later, in Example \ref{example48}, we will find the corresponding triple of $\tilde{Q}$, and we will show that $\tilde{Q}$ is the corresponding Weyl function.

2.2 Most of the statements in the first part of the following theorem about the Weyl function $Q$ are already known, as cited. We added a proof of regularity of $Q$ in order to obtain a characterization. Part (b) is more interesting. In part (b) we start from a generalized Nevanlinna function $Q$ and under the condition of regularity of $Q$ we prove the existence of a simple closed operator $S$ so that $Q$ becomes a Weyl function of $S$. Part (b) is a generalization of the converse part of \cite[Theorem 1]{DM1}.

\begin{thm}\label{theorem22} 

\begin{enumerate}[(a)]
\item Let $S$, $\lbrace 0 \rbrace \subseteq S \subsetneq A$, be a simple closed symmetric operator in a Pontryagin space $\mathcal{K}$ of index $\kappa$. Let $A^{+}=A$, $\rho(A) \neq \emptyset $, let $\Pi =(\mathcal{H},\Gamma_{0}, \Gamma_{1})$ be an ordinary boundary triple for $S^{+}$ ($A=\ker \Gamma_{0}$), and let $Q(z)$ be the Weyl function of $A$ corresponding to $\Pi$. Assume that $Q(w)$ is invertible for at least one point $w \in \mathcal{D}(Q)$.

Then $Q\in N_{\kappa }\left( \mathcal{H} \right)$, $Q$ is a regular and strict function uniquely determined by the relation $A$ in the minimal representation of the form (\ref{eq12}).

\item Conversely, let $Q\in N_{\kappa }\left( \mathcal{H} \right)$ be a regular and strict function given by a minimal representation (\ref{eq12}) with a representing relation $A$. 

Then there exists a unique closed simple linear operator $S$, $\lbrace 0\rbrace \subseteq S \subsetneq A \subsetneq S^{+} $ and there exists an ordinary boundary triple $\Pi =(\mathcal{H},\Gamma_{0}, \Gamma_{1})$ for $S^{+}$ such that $A=\ker \Gamma_{0}$. The function $Q(z)$ is the Weyl function of $A$ corresponding to $\Pi$.

\item In this case, the following hold:

\begin{enumerate}[(i)]
\item The representing relation $\hat{A}$ of $\hat{Q}:=-Q^{-1}$ satisfies $\hat{A} = \ker \Gamma_{1}$.
\item $A$ and $\hat{A}$ are transfersal extensions of $S:=A\cap \hat{A}$.
\end{enumerate}
\end{enumerate}
\end{thm}

\begin{proof} (a) The assumptions are appropriate. Namely, the existence of the boundary triple $\Pi =(\mathcal{H},\Gamma_{0}, \Gamma_{1})$, with $A:=\ker \Gamma_{0}$, has been proven in \cite[Proposition 2.2 (2)]{D1}. The existence of the corresponding (well defined) Weyl function with bounded values $Q(z)$ has been proven in \cite[p. 4427]{D1}. 

According to the terminology of \cite[p. 619]{BeLu}, the assumption that the closed linear relation $S$ is \textit{simple} means 
\begin{equation}
\label{eq21}
\mathcal{K}=c.l.s.\left\{\mathcal{R}_{z }(S^{+}):z\in \rho \left( A \right) \right\}. 
\end{equation}

The relationship between one-to-one operators $\gamma_z \in [\mathcal{H}, \mathcal{R}_{z }], z\in \rho(A)$, of the $\gamma$-field $\gamma$ and the Weyl function has been established by \cite[(2.13)]{D1}
\begin{equation}
\label{eq21a}
\frac{Q\left( z \right)-{Q^{*}\left( w \right)}}{z-\bar{w}}=\gamma_{w}^{+}\gamma_z, \thinspace \forall w, z \in \rho(A),  
\end{equation}
where, according to \cite[(2.6)]{D1}, $\gamma$-filed satisfies
\begin{equation}
\label{eq21b}
\gamma_{z}=\left( I+\left( z-w \right)\left( A-z \right)^{-1} \right)\gamma_{w}.
\end{equation}
For all $h, k \in \mathcal{H}$,
\[
\left(\frac{Q\left( z \right)-{Q^{*}\left( w \right)}}{z-\bar{w}}h,k\right)=\left(\gamma_{w}^{+}\gamma_z(h),k\right)=\left[\gamma_z(h), \gamma_{w}(k) \right]=\left[f,g\right], f \in \mathcal{R}_{z }, g\in \mathcal{R}_{w}.
\]
Because $\left(\mathcal{K}, [.,.] \right)$ given by (\ref{eq21}) is a Pontryagin space with a negative index $\kappa$, we conclude that $Q$ has $\kappa$ negative squares. Because $Q(z)$ are bounded operators, $Q\in N_{\kappa }\left( \mathcal{H} \right)$ holds. 

Let us note that the corresponding claim for Weyl families and generalized Nevanlinna families has been proven in \cite[Theorem 4.8]{BDHS}.  

From (\ref{eq21a}) and (\ref{eq21b}) it follows that
\begin{equation}
\label{eq21c}
Q\left( z \right)=Q(\bar{w})+(z-\bar{w})\gamma_{w}^{+}\left( I+\left( z-w \right)\left( A-z \right)^{-1} \right)\gamma_{w},z\in \rho\left( A \right).
\end{equation}

Because $\gamma_{z}(\mathcal{H})=\mathcal{R}_{z }$, according to (\ref{eq21}) and (\ref{eq21b}), the minimality condition (\ref{eq14}) is fulfilled with $ A = \ker \Gamma_{0}$ and with $\gamma$-field (\ref{eq21b}). Then, according to Theorem \ref{theorem12}, the state space $\mathcal{K}$, the representing relation $A$, the $\gamma$-field and the function $Q$ given by (\ref{eq21c}) are uniquely determined (up to unitary equivalence). 

By the definition of a $\gamma$-field, the operators $\gamma_{z}:\mathcal{H} \rightarrow \mathcal{R}_{z }$ are one-to-one for all $z \in \mathcal{D}(Q)$. Then, according to Proposition \ref{proposition22} (i), the function $Q(z)$ is strict. 

Let us prove that the function $Q$ is regular. According to our assumptions, there exists at least one point $\bar{w} \in \mathcal{D}(Q)$ such that $\hat{Q}(\bar{w}):=-Q(\bar{w})^{-1}$ is an operator. Because of the symmetry of the function $Q$, $Q(w)^{-1}$ is also an operator. According to definition (\ref{eq126}) of the Weyl function, it is obvious that $\mathcal{D}(\hat{Q}(z))= \ran\, \Gamma_{1} =\mathcal{H}, \forall z\in \mathcal{D}(Q)$. Therefore $\left( -Q(w)^{-1}\right)^{*}= \left( -Q(w)^{*}\right)^{-1}=\left( -Q(\bar{w})\right)^{-1}$ is an operator. This further means that $\hat{Q}(w)$ is a closed operator. It is also defined on entire $\mathcal{H}$, i.e., $\hat{Q}(w)$ is bounded operator. This proves that $Q(w)$ is boundedly invertible operator. By definition $Q$ is a regular function. This completes the proof of (a).
\\

(b) The assumption that $Q\in N_{\kappa }\left( \mathcal{H} \right)$ is a regular function with the representing relation $A$ in the minimal representation (\ref{eq12}) includes that (\ref{eq14}) and (\ref{eq16}) hold, and $\rho (A)\ne \emptyset$. According to \cite [Proposition 2.1] {Lu1}, the inverse $\hat{Q}=-Q^{-1} \in N_{\kappa }\left( \mathcal{H} \right)$ admits the representation 
\begin{equation}
\label{eq22}
\hat{Q}\left( z \right)=\hat{Q}\left( \overline{w} \right)+\left( 
z-\overline{w} \right)\hat{\Gamma }^{+}_{w}\left( I+\left( z-w \right)(\hat{A}-z)^{-1} \right)\hat{\Gamma}_{w},
\end{equation}
where $ w\in \rho \left( A \right)\cap \rho ( \hat{A} )$ is an arbitrarily selected point of reference,
\begin{equation}
\label{eq23}
\hat{\Gamma}_{w}:=-\Gamma_{w} {Q(w)}^{-1}, 
\end{equation}
and 
\begin{equation}
\label{eq24}
( \hat{A}-z )^{-1}=\left( A-z \right)^{-1}-\Gamma_{z}{Q\left( z 
\right)}^{-1}\Gamma_{\bar{z}}^{+},\, \, \forall z\in \rho \left( A 
\right)\cap \rho ( \hat{A})
\end{equation}
holds.

According to Proposition \ref{proposition22} (i), the assumption that $Q\in N_{\kappa }\left( \mathcal{H} \right)$ is a strict function means that operators $\Gamma_{z}, z\in \mathcal{D}(Q)$, in representation (\ref{eq12}) are one-to-one. 

We need to prove that there exists a closed symmetric relation $S$, a boundary triple $\Pi=\left( \mathcal{H}, \Gamma_{0}, \Gamma_{1}\right) $ and a corresponding Weyl function $M(z)=Q(z)$. 

We define the closed symmetric relation $S$ by
\begin{equation}
\label{eq24a}
S:=A\cap \hat{A}.
\end{equation}
Because representations (\ref{eq12}) and (\ref{eq22}) are uniquely determined, the linear relation $S$ is also uniquely determined. This also means that the self-adjoint relation $A$ is an extension of $S$. 

The linear relation $S$ defined by (\ref{eq24a}) has equal (finite or infinite) defect numbers in the separable Pontryagin space $\mathcal{K}$ because it has a self-adjoint extension $A$ within $\mathcal{K}$. Let us denote that defect number by $d(S)$. We already observed that $\Gamma_{z}: \mathcal{H}\rightarrow \Gamma_{z}(\mathcal{H}), \, z\in \rho(A)$, are one-to-one operators. Therefore, $ \dim \mathcal{H}=d(S)$. 

We can here apply \cite[Proposition 2.2]{D1}. Therefore, there exists a boundary triple $\tilde{\Pi} =(\mathcal{\tilde{H}},\Gamma_{0}, \Gamma_{1})$ for $S^{+}$ such that $A=\ker \Gamma_{0}$, with a $\gamma$-field $\gamma_{z}, z \in \rho(A)$, that satisfies (\ref{eq21b}). 

According to \cite[Proposition 2.2 (3)]{D1}, $\gamma_{z}:\tilde{\mathcal{H}}\rightarrow \mathcal{R}_{z}=\gamma_{z}( \tilde{\mathcal{H}}) , \forall z\in \rho(A)$, is a one-to-one operator. Recall that $\gamma_{z}$ and $\tilde{\mathcal{H}} $ were introduced so that $\dim\, (\tilde{\mathcal{H}})=d(S)$ holds. This means $\dim\, (\tilde{\mathcal{H}})=\dim\, \mathcal{H}=d(S)$. Therefore, we can  consider $\mathcal{H}=\tilde{\mathcal{H}}$, hence $\tilde{\Pi} =(\mathcal{H},\Gamma_{0}, \Gamma_{1})$.

Let $M(z)$ be the Weyl function corresponding to $\tilde{\Pi} =(\mathcal{H},\Gamma_{0}, \Gamma_{1})$.
Then $M(z)$ and $\gamma (z)$ satisfy \cite[(2.13)]{D1}. According to \cite[Remark 2.2]{D1}, the operator valued function $M(z)$ is a $Q$-function of $S$ represented by $A=\ker\,\Gamma_{0}$ in some Pontryagin space $\mathcal{\tilde{K}}$. (For a definition of the $Q$-function of $S$ see e.g. \cite{KL3}.) The minimal Pontryagin space of the $Q$-function $M(z)$ is given by means of $\gamma_{z}(\mathcal{H})=\mathcal{R}_{z }(S^{+})$, which is
\begin{equation}
\label{eq25a}
\tilde{\mathcal{K}}:=c.l.s.\left\{\mathcal{R}_{z }(S^{+}):z\in \rho \left( A \right) \right\} \subseteq \mathcal{K}.
\end{equation}

According to \cite[(2.13)]{D1} and (\ref{eq21b})
\begin{equation}
\label{eq214}
M\left( z \right)={M(w)}^{\ast }+(z-\bar{w})\gamma_{w}^{+}\left( I+\left( z-w \right)\left( A-z \right)^{-1} \right)\gamma_{w}, z\in \rho\left( A \right).
\end{equation}

Let us now use the so called $\epsilon_{z}$-model, see \cite{Lu3, KL2}. According to that model, we can identify the building blocks of $\tilde{\mathcal{K}}$ with $\gamma_{z}(h)\, (h\in  \mathcal{H},\, z\in \rho(A))$, and the building blocks of $\mathcal{K}$ with $\Gamma_{z}(h), (h\in \mathcal{H}, \, z\in \rho(A))$. Therefore, we can define one-to-one operator $U: \tilde{\mathcal{K}}\rightarrow\mathcal{K}$ by 
\[
U\left( \gamma_{z}(h)\right) =\Gamma_{z}(h),\, \forall h\in \mathcal{H}, \, \forall z\in \rho(A),
\]
and we can set 
\[
\left[\gamma_{z}(h), \gamma_{w}(k) \right] = \left[\Gamma_{z}(h), \Gamma_{w}(k) \right],\, \forall h, k\in \mathcal{H}, \, \forall z,w\in \rho(A). 
\]
Obviously, the operator $U$ is a unitary operator. Therefore, the spaces $\tilde{\mathcal{K}} $ and $\mathcal{K} $ are unitarily equivalent. This, together with $\mathcal{H}=\tilde{\mathcal{H}}$, means that the representations (\ref{eq12}) and (\ref{eq214}), both represented by the same relation $A$, are unitarily equivalent. In other words, we can consider $Q=M$.

According to (\ref{eq25a}), by definition $S$ is a simple relation with respect to $\tilde{\mathcal{K}}=\mathcal{K}$. We know that a simple linear relation $S$ is an operator. 

(c) (i) According to \cite [(2.3)]{D1}, there exists a bijective correspondence between proper extensions $\tilde{S}\in Ext\, S$ and closed sub-spaces $\theta $ in $\mathcal{H}\times \mathcal{H}$ defined by
\begin{equation}
\label{eq26}
S_{\theta } \in \, Ext\, S\, \Leftrightarrow \, \theta 
:=\Gamma S_{\theta }=\left\{ \left( {\begin{array}{*{20}c}
\Gamma_{0}\hat{f}\\
\Gamma_{1}\hat{f}\\
\end{array} } \right): \hat{f}\in S_{\theta } \right\}\in \tilde{\mathcal{C}}(\mathcal{H}).
\end{equation}
Then the Krein (a.k.a. Krein-Naimark) formula 
\begin{equation}
\label{eq28}
\left( S_{\theta }-z \right)^{-1}=\left( A-z \right)^{-1}+\Gamma_{z}\left( 
\theta -Q\left( z \right) \right)^{-1}\Gamma_{\bar{z}}^{+}
\end{equation}
holds. Let us set $S_{\theta }:=\hat{A}$, where $\hat{A}$ is the linear relation that represents the inverse function $\hat{Q}$ in representation (\ref{eq22}). Then according to (\ref{eq24}), the pair:  $S_{\theta }=\hat{A},  \theta =O_{\mathcal{H}}$ (a zero function on $\mathcal{H}$), satisfies (\ref{eq28}). Because the correspondence defined by (\ref{eq26}) is a bijection, it follows 
\begin{equation}
\label{eq210}
\theta =\Gamma \hat{A}=\left\{ \left( {\begin{array}{*{20}c}
\Gamma_{0}\hat{f}\\
0\\
\end{array} } \right) : \hat{f}\in \hat{A} \right\}.
\end{equation}
Therefore, $\hat{A}=\ker \Gamma_{1}=:S_{1}$. This proves (ii). 

(ii) $S:=A\cap \hat{A}$ has been defined in (b). It suffices to prove $S^{+} \subseteq \ker \Gamma_{0} \hat{+} \ker \Gamma_{1}$. 

Assume $\hat{k} \in S^{+}$ and $\hat{h}=\Gamma \hat{k}$. Then, because $\Gamma$ is surjective, we have
\[
\left( {\begin{array}{*{20}c}
h\\
h'\\
\end{array} } \right)=\left( {\begin{array}{*{20}c}
0\\
h'\\
\end{array} } \right)+\left( {\begin{array}{*{20}c}
h\\
0\\
\end{array} } \right)=\Gamma \hat{t}+\Gamma \hat{r}, \, \hat{t} \in \ker\, \Gamma _{0}, \hat{r} \in \ker\, \Gamma _{1}.
\]
Hence, $ \hat{s}:=\hat{k}-\hat{t}-\hat{r} \in S \subseteq ker\, \Gamma _{0}$, i.e.
$\hat{k}:=\left( \hat{s}+ \hat{t}\right)+\hat{r} =: \hat{u}+ \hat{r} \in \ker\, \Gamma _{0} \hat{+} \ker\, \Gamma _{1}$. This proves $S^{+} \subseteq \ker \Gamma_{0} \hat{+} \ker \Gamma_{1}$.
\end{proof}

\begin{cor}\label{corollary24}. Let $\mathcal{K}$ be a Pontryagin space of negative index $\kappa$ and let $M(z)$ be the Weyl function associated with the ordinary boundary triple $\Pi =(\mathcal{H},\Gamma_{0}, \Gamma_{1})$. If $\hat{M}:=-M^{-1}$ exists then relations $S_{i}:=\ker \Gamma_{i}, i=1,2$, satisfy 
\begin{equation}
\label{eq216}
( S_{1}-z )^{-1}=\left( S_{0}-z \right)^{-1}+\hat{\gamma}_{z} \gamma_{\bar{z}}^{+}, z\in \rho \left(S_{0} \right)\cap \rho ( S_{1} ),
\end{equation}
where $\gamma_{z}$ and $\hat{\gamma}_{z}$ are $\gamma$-fields associated with $S_{0}$ and $S_{1}$, respectively.
\end{cor}
\begin{proof} By definition of the Weyl function, the operator $\Gamma_{1}$ is for $\hat{M}$ what $\Gamma_{0}$ is for $M$. According to Theorem \ref{theorem22} (c), $\hat{A}=S_{1}$. Therefore, we can substitute $S_{0}$ and $S_{1}$ for $A$ and $\hat{A}$ into (\ref{eq24}), respectively. Hence, we can rewrite (\ref{eq23}) with $w=z$, $\Gamma_{w}=\gamma_{z} $, $\hat{\Gamma}_{w}=\hat{\gamma}_{z}$ and substitute (\ref{eq23}) into (\ref{eq24}) to obtain (\ref{eq216}). 
\end{proof}

\textbf{2.3.} Identity (\ref{eq216}) gives us a relationship between resolvents of $A=\ker \Gamma_{0}$ and $\hat{A}:=\ker \Gamma_{1}$ when $S:=A \cap\hat A$ and $A$ is the representing relation of the Weyl function $Q$, i.e. of the regular and strict generalized Nevanlinna function $Q$. In the following proposition, we will establish a direct relationship between any two closed linear relations $A$ and $B$ that satisfy $\rho(A) \cap \rho(B) \neq \emptyset$. Then we will apply it to the representing relations $A$ and $\hat A$ of $Q$ and $\hat{Q}$, respectively.

Recall, for the \textit{defect subspace} of a linear relation $T$ we use the notation 
\[
\hat{\mathcal{R}}_{z}\left( T \right)=\left\{ \left( {\begin{array}{*{20}c}
t\\
zt\\
\end{array} } \right)\in T \right\}.
\]
\begin{propo}\label{proposition26}. Let $A$ and $B$ be linear relations in a Krein space $\mathcal{K}$, let $B$ be a closed relation, and $\rho \left( A \right) \cap \rho \left( B \right)\ne \emptyset$. Then 
\begin{equation}
\label{eq218}
A \subseteq B\dotplus \hat{\mathcal{R}}_{z}\left( A \hat{+} B \right), \forall \, z\in \rho \left( A \right) \cap \rho \left( B \right).
\end{equation}
Equality holds if and only if $A=B$.
\end{propo}

\begin{proof} For $z\in \rho \left( A \right) \cap \rho \left( B \right)$ and for every 
$\left( {\begin{array}{*{20}c}
f\\
f^{'}\, \\
\end{array} } \right)\in A$ we have $ \left( {\begin{array}{*{20}c}
f\\
f^{'}-zf\, \\
\end{array} } \right)\in A-z $.
Because $z\in \rho \left( B \right)$, and $B$ is closed, there exists $\left( {\begin{array}{*{20}c}
g\\
g^{'}\\
\end{array} } \right) \in B$ such that 
\[
f^{'}-zf\, =g^{'}-zg\Rightarrow f^{'}-g^{'}=z\left( f-g \right)
\]
holds. Therefore 
\[
\left( {\begin{array}{*{20}c}
f\\
f^{'}\\
\end{array} } \right)-\, \left( {\begin{array}{*{20}c}
g\\
g^{'}\\
\end{array} } \right)=\left( {\begin{array}{*{20}c}
f-g\\
f^{'}-g^{'}\\
\end{array} } \right)=\left( {\begin{array}{*{20}c}
f-g\\
z\left( f-g \right)\\
\end{array} } \right) \in \hat{\mathcal{R}}_{z}\left( A\hat{+}B \right).
\]
Thus
\begin{equation}
\label{eq220}
\left( {\begin{array}{*{20}c}
f\\
f^{'}\\
\end{array} } \right)= \left( {\begin{array}{*{20}c}
g\\
g^{'}\\
\end{array} } \right)+\left( {\begin{array}{*{20}c}
f-g\\
z\left( f-g \right)\\
\end{array} } \right).
\end{equation}

The sum (\ref{eq220}) is direct because $0\ne \left( {\begin{array}{*{20}c}
t\\
zt\\
\end{array} } \right)\in B \cap \hat{\mathcal{R}}_{z}\left(A\hat{+}B 
\right)\Rightarrow z\in \sigma_{p}\left( B\right)$, which contradicts the assumption $z\in  \rho \left( B \right)$. This proves (\ref{eq218}). 

To prove the second claim, let us assume $A = B \dotplus \hat{\mathcal{R}}_{z}\left( A \hat{+} B \right), z\in \rho \left( A \right) \cap \rho \left( B \right)$. Then for $S:=A\cap B$ we have
\[
S=B \subseteq A \Rightarrow A \hat{+} B = A \Rightarrow \hat{\mathcal{R}}_{z}\left( A \hat{+} B \right)=\emptyset \Rightarrow A=B.
\] 
The converse implication follows from $\hat{\mathcal{R}}_{z}\left( B \right)=\lbrace 0 \rbrace $. 
\end{proof}

\begin{cor}\label{corollary28}. Let $Q\in N_{\kappa }(\mathcal{H})$ be a regular strict function and let $A$ and $\hat{A}$ be the representing relations of $Q$, and $\hat{Q}:=-Q^{-1}$, respectively. For $S=A\cap \hat{A}$, 
\[
A \subseteq \hat{A}\dot{+}\hat{\mathcal R}_{z}\left(S^{+} \right), \forall  z\in \rho( A ) \cap \rho( \hat{A} ).
\]
holds. Equality holds if and only if $A=\hat{A}$.
\end{cor}

\textbf{Proof.} The regularity of $Q$ implies $\rho \left( A \right)\cap \rho ( \hat{A} )\ne \emptyset $. According to Theorem \ref{theorem22} (c)(ii), we can substitute $S^{+}$ for $A \hat{+} \hat{A}$. Then both claims follow from Proposition \ref{proposition26}.\hfill $\square$ 

Obviously, the relations $A$ and $\hat{A}$ can exchange places in the above corollary.

\section{Weyl function $Q \in N_{\kappa}(\mathcal{H})$ with  boundedly invertible $Q^{'}(\infty)$}\label{s6} 
3.1 A significant part of this paper is about the class of functions $Q\in N_{\kappa }(\mathcal{H})$ that are holomorphic at $\infty $, i.e. the functions $Q$ for which there exists $Q^{'}\left( \infty \right):={\lim \limits_{z \to \infty}{zQ(z)}}$.

\begin{lem}\label{lemma14} \cite[Lemma 3]{B} A function $Q \in N_{\kappa }(\mathcal{H})$ is 
holomorphic at $\infty $ if and only if $Q(z)$ has a representation
\begin{equation}
\label{eq110}
Q\left( z \right)=\tilde{\Gamma}^{+}\left( A-z \right)^{-1}\tilde{\Gamma} ,z\in \rho \left( A \right),
\end{equation}
with a bounded operator $A$. In this case 
\begin{equation}
\label{eq114}
Q^{'}\left( \infty \right):={\lim \limits_{z \to \infty}{zQ(z)}}=-\tilde{\Gamma}^{+}\tilde{\Gamma},
\end{equation}
where the limit denotes convergence in the Banach space of bounded operators $\mathcal{L}(\mathcal{H})$.
\end{lem}

Recall, see \cite[Proposition 1]{B}, that the operator $\tilde{\Gamma} $ used in (\ref{eq110}) can be expressed as 
\begin{equation}
\label{eq116}
\tilde{\Gamma} =\left( A-z \right)\Gamma_{z},  \forall z \in \rho \left( A \right ).
\end{equation}
Then the representation (\ref{eq110}) is minimal, if and only if 
\[
\mathcal{K}=c.l.s.\left\{ \left( A-z \right)^{-1}\tilde{\Gamma}h :z\in \rho \left( A \right),h\in \mathcal{H} \right\}.
\]
The decomposition of the function $Q \in N_{\kappa }(\mathcal{H})$ in \cite[Remark 1]{B} shows us the important role representations of the form (\ref{eq110}) play in research of the function $Q\in N_{\kappa}(\mathcal{H})$. 

The following lemma from \cite{B} will be frequently needed in this paper.

\begin{lem}\label{lemma16} \cite[Lemma 4]{B} Let $\tilde{\Gamma}:\mathcal{H}\to \mathcal{K}$ be a bounded operator and let $\tilde{\Gamma}^{+}:\mathcal{K}\to \mathcal{H}$ be its adjoint operator. Assume also that $\tilde{\Gamma}^{+}\tilde{\Gamma} $ is a boundedly invertible operator in the Hilbert space $\left( \mathcal{H},\left( .,. \right) 
\right)$. Then for the operator 
\begin{equation}
\label{eq118}
P:=\tilde{\Gamma} \left( \tilde{\Gamma}^{+}\tilde{\Gamma} \right)^{-1}\tilde{\Gamma}^{+}
\end{equation}
the following statements hold:
\begin{enumerate}[(i)]
\item P is an orthogonal projection in the Pontryagin space $(\mathcal{K},\, \left[ .,. \right]).$
\item The scalar product $\left[ .,.\right] $ does not degenerate on $P\mathcal{K}=\tilde{\Gamma} \mathcal{H} $ and therefore it does not degenerate on ${\tilde{\Gamma} \left( \mathcal{H} \right)\, }^{[\bot ]}=\ker \tilde{\Gamma}^{+}. $
\item $\ker \tilde{\Gamma}^{+}=\left( I-P \right)\mathcal{K}.$ 
\item The Pontryagin space $\mathcal{K}$ can be decomposed as a direct orthogonal sum of Pontryagin spaces i.e. 
\end{enumerate}
\begin{equation}
\label{eq120}
\mathcal{K}=\left( I-P \right)\mathcal{K}[+]P\mathcal{K}.
\end{equation}
\end{lem}
3.2 Let 
\[
\mathcal{K}:=\mathcal{K}_{1}\left[ + \right]\mathcal{K}_{2}
\]
be a Pontryagin space with nontrivial Pontryagin subspaces $\mathcal{K}_{l},\, l=1,\, 
2$, and let $E_{l}:\mathcal{K}\to \mathcal{K}_{l},\, l=1,\, 2$, be orthogonal projections. Let $T$ be a linear relation in $\mathcal{K}=\mathcal{K}_{1}\left[ + \right]\mathcal{K}_{2}$. If for any projection $E_{i},\, i=1, 2$, $E_{i}\left( D(T)\right)\subseteq D(T)$ holds, then according to \cite[Lemma 2.2]{B1} the following four linear relations can be defined 
\[
T_{i}^{j}:=\left\{ \left( {\begin{array}{*{20}c}
k_{i}\\
k_{i}^{j}\\
\end{array} } \right):k_{i}\in D\left( T \right)\cap \mathcal{K}_{i},\, k_{i}^{j}\in 
E_{j}T\left( k_{i} \right) \right\}\subseteq \mathcal{K}_{i}\times \mathcal{K}_{j},\, i, j=1,\, 2.
\]
In this notation the subscript ``$i$'' is associated with the domain subspace $\mathcal{K}_{i}$, the superscript ``$j$'' is associated with the range subspace $\mathcal{K}_{j}$. For example $\left( {\begin{array}{*{20}c}
k_{1}\\
k_{1}^{2}\\
\end{array} } \right)\in T_{1}^{2}$. We will use ``$\left[+\right]$'' 
to denote adjoint relations of $T_{i}^{j}$. Therefore 
\[
T_{1}^{2}\subseteq \mathcal{K}_{1} \times \mathcal{K}_{2}\Rightarrow { T_{1}^{2}}^{\left[+\right]}\subseteq \mathcal{K}_{2} \times \mathcal{K}_{1}.
\]
Hence, for the linear relation $T$ and decomposition $\mathcal{K}:=\mathcal{K}_{1}\left[ + \right]\mathcal{K}_{2}$, we can assign the following \textit{relation matrix} 
\[
\left( {\begin{array}{*{20}c}
T_{1}^{1} & T_{2}^{1}\\
T_{1}^{2} & T_{2}^{2}\\
\end{array} } \right).
\]
We obtain 
\[
T=\left( T_{1}^{1}+T_{1}^{2} \right)\hat{+}\left( T_{2}^{1}+T_{2}^{2} 
\right).
\] 
\begin{lem}\label{lemma32}. Let $Q \in N_{\kappa }\left( \mathcal{H} \right)$ satisfy conditions of Lemma \ref{lemma14}. Then 
\begin{equation}
\label{eq312}
B:=A_{\vert (I-P)\mathcal{K}}\dot{+}\left( \left\{ 0 \right\}\times P \mathcal{K} \right)\subseteq (I-P)\mathcal{K} \times \mathcal{K}
\end{equation}
holds, where projection $P$ is defined by (\ref{eq118}).
Then 
\begin{equation}
\label{eq314}
z\in \rho \left( A \right)\cap \rho ( \hat{A} )\Rightarrow \mathcal{K}\subseteq \left( B-z \right)\left( I-P \right)\mathcal{K},
\end{equation}
and 
\[
z\in \rho \left( A \right)\cap \rho ( \hat{A} )\Rightarrow z\in \rho \left( B \right).
\]
\end{lem}
\begin{proof} Assume $z\in \rho \left( A \right)\cap \rho ( \hat{A} )$. Then, according to (\ref{eq22}) and \cite[Theorem 3]{B}, $z \in \rho ( \hat{A} )$ if and only if $ z \in \rho (\tilde{A} )$, where 
\[
\tilde{A}:= (I-P)A_{\vert (I-P)\mathcal{K}}.
\]
Therefore, for any $ f=\left( I-P \right)f+Pf\in \mathcal{K} $ there exists $g\in (I-P)\mathcal{K}$, such that 
\[
\left( I-P \right)f=\left( {\left( I-P \right)A}_{\vert (I-P)\mathcal{K}}-z\left( I-P \right) \right)g.
\]
Also, there exists $k\in \mathcal{K}$ such that 
\[
Pk=Pf-PA_{\vert (I-P)\mathcal{K}}g\Rightarrow Pf=PA_{\vert (I-P)\mathcal{K}}g+Pk
\]
holds. We will also use the identity: ${\left( I-P \right)A}_{\vert (I-P)\mathcal{K}}+PA_{\vert (I-P)\mathcal{K}}=A_{\vert (I-P)\mathcal{K}}$.

Now we have, 
\[
f=\left( I-P \right)f+Pf
\]
\[
=\left( {\left( I-P \right)A}_{\vert (I-P)\mathcal{K}}-z\left( I-P \right) \right)g+PA_{\vert (I-P)\mathcal{K}}g+Pk
\]
\[
=\left( A_{\vert (I-P)\mathcal{K}}-z\left( I-P \right) \right)g+Pk \in \left( B-z(I-P) \right)g \in \left( B-z \right)\left( I-P \right)\mathcal{K}.
\]
This proves (\ref{eq314}).

Let us prove that for $ z\in \rho \left( A \right)\cap \rho (\hat{A})$ and $f\in (I-P)\mathcal K$  
\[
(B-z)f=0 \Rightarrow f=0
\] 
holds. Indeed, we already mentioned that $ z\in \rho \left( A \right)\cap \rho (\hat{A}) \Rightarrow  z\in \rho \left( (I-P)A_{\vert (I-P)\mathcal{K}} \right)$. Now we have for $f\in (I-P)\mathcal{K}$:
\[
0=(B-z)f=\left( A_{\vert (I-P)\mathcal{K}}\dot{+}\left( \left\{ 0 \right\}\times P \mathcal{K} \right) - z \right) f \Rightarrow
\]
\[
\Rightarrow \left( (I-P)A_{\vert (I-P)\mathcal{K}} -z \right)f=0. 
\]
Because $z\in \rho \left( (I-P)A_{\vert (I-P)\mathcal{K}} \right)$, it follows that $f=0$. This further means that $(B-z)^{-1}$ is an operator. Relation $B$ is closed as a sum of a bounded and closed relation. Then, because of (\ref{eq314}) the closed operator $(B-z)^{-1}$ is bounded. This proves $ z\in \rho \left( B \right)$.
\end{proof}
Now we can prove the following lemma. 
\begin{lem}\label{lemma34} Let $Q\in N_{\kappa }\left( \mathcal{H} \right)$ satisfy conditions of Lemma \ref{lemma14}. Then the representing relation $\hat{A}$ of $\hat{Q}:=-Q^{-1}$ satisfies
\begin{equation}
\label{eq316}
\hat{A}=A_{\vert (I-P)\mathcal{K}}\dot{+}\hat{A}_{\infty },
\end{equation}
where 
\[
\hat{A}_{\infty }= \left\{ 0 \right\}\times P \mathcal{K}.
\]
\end{lem}
\begin{proof} Because $\tilde{\Gamma}^{+}\tilde{\Gamma} $ is boundedly invertible, according to Lemma \ref{lemma16} the scalar product $[.,.]$ does not degenerate on the subspace $P\left( \mathcal{K} \right)=\tilde{\Gamma} (\mathcal{H})$. According to \cite[Theorem 3]{B}, there exists 
$\hat{Q}\left( z \right):=-{Q\left( z \right)}^{-1}, z\in \rho\left( A \right)\cap (\hat{A})$. Let $\hat{Q}$ be represented by a self-adjoint linear relation $\hat{A}$ in representation (\ref{eq22}). Then $\hat A$ satisfies (\ref{eq24}). 

Let us now observe the linear relation $B$ given by (\ref{eq312}), and let us find the resolvent $\left( B-z \right)^{-1}$, which exists according to Lemma \ref{lemma32}. Let us select a point $z\in \rho (B)$ and a vector 
\[
f\in \mathcal{K}=\left( B-zI \right)\left( I-P \right)\mathcal{K},
\]
and let us find $ \left( B-z \right)^{-1}f $.

According to Lemma \ref{lemma32} there exists an element $g:=\left( B-z \right)^{-1}f \in \left( I-P \right)\mathcal{K}$. According to definition (\ref{eq312}) of $B$ and $P\mathcal{K}=\tilde{\Gamma} \mathcal{H} $,
\[
\left\{ g, f+zg \right\}\in A_{\vert (I-P)\mathcal{K}}\dot{+}\left( \left\{ 0 \right\}\times \tilde{\Gamma} \mathcal{H}\right)
\]
holds. This means that for some $h\in \mathcal{H}$ 
\[
f+zg=Ag + \tilde{\Gamma} h
\]
holds. Then we have
\[
Ag-zg=f-\tilde{\Gamma} h.
\]
Hence,
\[
g=\left( A-z \right)^{-1}f-\left( A-z \right)^{-1}\tilde{\Gamma} h.
\]
Because $\tilde{\Gamma}^{+}(I-P)=0$, we have
\[
0=\tilde{\Gamma}^{+}g= \tilde{\Gamma}^{+}\left( A-z \right)^{-1}f-\tilde{\Gamma}^{+}\left( A-z \right)^{-1}\tilde{\Gamma} h=\tilde{\Gamma}^{+}\left( A-z \right)^{-1}f-Q\left( z \right)h.
\]
According to (\ref{eq116}),
\[
\Gamma_{\bar{z}}^{+}f=\tilde{\Gamma}^{+}\left( A-z \right)^{-1}f.
\]
Therefore,
\[
h={Q\left( z \right)}^{-1}\Gamma_{\bar{z}}^{+}f.
\]
This and (\ref{eq116}) gives 
\[
\left( B-z \right)^{-1}f=g=\left( A-z \right)^{-1}f-\Gamma_{z}h =\left( A-z \right)^{-1}f-\Gamma_{z}{Q\left( z 
\right)}^{-1}\Gamma_{\bar{z}}^{+}f,
\]
which proves that formula (\ref{eq24}) holds for the linear relation $B \subseteq \left(I-P \right) \mathcal{K}\times \mathcal{K}$ defined by (\ref{eq312}). Therefore, $\left( B-z \right)^{-1}=( \hat{A}-z )^{-1}$, and 
\[
\hat{A}=B=A_{\vert (I-P)\mathcal{K}}\dot{+}\left( \left\{ 0 \right\}\times P\mathcal{K} \right).
\]
Because $A$ is a single valued, the sum is direct, and $\hat{A}_{\infty }=\left( \left\{ 0 \right\}\times P \mathcal{K} \right)$, i.e. representation (\ref{eq316}) of $\hat{A}$ holds.
\end{proof}

Note, identity (\ref{eq316}) derived here by means of the operator valued function $Q\in N_{\kappa }\left( \mathcal{H} \right)$ corresponds to identity \cite [(3.5)] {HL} which was derived for a scalar function $ q \in N_{\kappa } $. Also note that $\hat{A}_{\infty }= \left\{ 0 \right\}\times P \mathcal{K}$ holds according to \cite[Proposition 5]{B} too. 

\begin{thm}\label{theorem32} Let $Q\in N_{\kappa }\left( \mathcal{H} \right)$ be holomorphic at infinity with boundedly invertible $Q^{'}\left( \infty \right)$ and let $Q$ be minimally represented by (\ref{eq110}) 
\[
Q\left( z \right)=\tilde{\Gamma}^{+}\left( A-z \right)^{-1}\tilde{\Gamma} , z\in \rho\left( A \right),
\]
with a bounded operator $A$. Then, relative to decomposition (\ref{eq120}) 
\[
\mathcal{K}_{1}\left[ + \right]\mathcal{K}_{2}:=\left( I-P \right)\mathcal{K}\left[ + \right]P\mathcal{K},
\]
the following hold:
\begin{enumerate}[(i)]
\item $A=\left( {\begin{array}{*{20}c}
\tilde{A} & (I-P)A_{\vert P\mathcal{K}}\\
PA_{\vert (I-P)\mathcal{K}} & PA_{\vert P\mathcal{K}}\\
\end{array} } \right)$, where $ \tilde{A}=(I-P)A_{\vert (I-P)\mathcal{K}}$.  
\item $\hat{A}=A_{\vert I-P} \dot{+} (\lbrace 0 \rbrace \times P\mathcal{K})=\left( {\begin{array}{*{20}c}
\tilde{A} & 0\\
0 & \hat{A}_{\infty }\\
\end{array} } \right)$, 
\item $S=A_{\vert (I-P)\mathcal{K}}$, $\mathcal{R}=P \mathcal{K}$. $S$ is a symmetric, closed, bounded operator.
\item $ S^{+}=\left( {\begin{array}{*{20}c}
\tilde{A} & (I-P)A_{\vert P\mathcal{K}}\\
(I-P)\mathcal{K} \times P\mathcal{K} & P \mathcal{K} \times P\mathcal{K}\\
\end{array} } \right)$. 
\item $\mathcal{R}_{z}=\left\lbrace \left( \begin{array}{*{20}c}
-(\tilde{A}-z)^{-1}APx_{P}\\
x_{P}\\
\end{array} \right) : x_{P} \in P\mathcal{K} \right\rbrace$, $\mathcal{K}=c.l.s. \left \lbrace \mathcal{R}_{z}: z \in \rho(A) \right \rbrace$, i.e. $S$ is simple. 
\item If additionally, $\tilde{\Gamma}$ is a one-to-one operator, then $Q$ is the Weyl function associated with $(S,A)$ and
$S^{+}=A\hat{+}\hat{A}=A\dot{+}\mathcal{\hat{R}}$.
\end{enumerate}
\end{thm}
\begin{proof} 
(i) The relation matrix of the operator $A$, with respect to decomposition (\ref{eq120}), is obviously
\begin{equation}
\label{eq322}
A=\left( {\begin{array}{*{20}c}
\left( I-P \right)A_{\vert (I-P)\mathcal{K}} & \left( I-P \right)A_{\vert P\mathcal{K}}\\
PA_{\vert (I-P)\mathcal{K}} & PA_{\vert P\mathcal{K}}\\
\end{array} } \right)=A_{\vert (I-P)\mathcal{K}}\dot{+}A_{\vert P\mathcal{K}}.
\end{equation}

(ii) According to \cite[Theorem 3 (ii)]{B}, the function $Q$ is regular. Therefore, there exists the inverse function $\hat{Q}$ and the representing relation $\hat{A}$.  According to (\ref{eq316}), the condition $(I-P)D(\hat{A})\subseteq D(\hat{A})$ is satisfied. Hence, according to \cite[Lemma 2.2]{B1}, there exists a relation matrix of $\hat{A}$ relative to decomposition (\ref{eq120}). Let that relation matrix be
\[
\hat{A}=\left( {\begin{array}{*{20}c}
{\hat{A}}_{1}^{1} & {\hat{A}}_{2}^{1}\\
{\hat{A}}_{1}^{2} & {\hat{A}}_{2}^{2}\\
\end{array} } \right),
\]
where ${\hat{A}\mathrm{\, }}_{i}^{j}\subseteq \mathcal{K}_{i} \times \mathcal{K}_{j},\, i, j=1,2$. 
According to Lemma \ref{lemma34}
\begin{equation}
\label{eq323}
\hat{A}\left( 0 \right)=P\mathcal{K}.
\end{equation}
Therefore, $\hat A(0)$ is an ortho-complemented subspace of $\mathcal{K}$. According to \cite [Theorem 2.4] {S}, 
\begin{equation}
\label{eq324}
\hat{A}=\hat{A}_{s}\dot{[+]}\hat{A}_{\infty },
\end{equation}
where $\hat{A}_{s}$ is a self-adjoint densely defined operator in ${\hat{A}(0)}^{[\bot ]}=(I-P)\mathcal{K}$, $\ran\, \hat{A}_{s} \subseteq \left( I-P \right)\mathcal{K}$ and denotes direct orthogonal sum of sub-spaces.

For $g\in \left( I-P \right)\mathcal{K}$, from (\ref{eq316}) and (\ref{eq324}), it follows that
\[
\left( \left( I-P \right)A_{\vert (I-P)\mathcal{K}}[\dot{+}]PA_{\vert (I-P)\mathcal{K}}\right)g + P k_{0}=A_{s}g[\dot{+}]Pk
\]
for some $k_{0}, k \in \mathcal{K}$. Obviously:
\begin{equation}
\label{eq326}
A_{s}g=\left( I-P \right)A_{\vert (I-P)\mathcal{K}}g=\tilde{A}g.
\end{equation}

Obviously $\hat{A}_{\infty } \subseteq P \mathcal K \times P \mathcal{K} $ and $\hat{A}_{\infty }={\hat{A}}_{\infty}^{+}$. Hence, the relation matrix of $\hat{A}$ is 
\begin{equation}
\label{eq328}
\hat{A}=\left( {\begin{array}{*{20}c}
\tilde{A} & 0\\
0 & \hat{A}_{\infty }\\
\end{array} } \right),
\end{equation}. 

(iii) Let us now find $S=A \cap \hat{A}$. According (\ref{eq328}), we have
\[
\hat{A}=\left\lbrace \left( {\begin{array}{*{20}c}
x_{I-P}\\
\tilde{A}x_{I-P}+p\\
\end{array} }\right): x_{I-P} \in (I-P)\mathcal{K}, \, p \in P\mathcal{K} \right\rbrace.
\]
Since $\dom\, S=(I-P)\mathcal{K}$, elements of $A\cap S$ satisfy
\[
\left( {\begin{array}{*{20}c}
x_{I-P}\\
Ax_{I-P}\\
\end{array} }\right)= \left(\begin{array}{*{20}c}
x_{I-P}\\
\tilde{A}x_{I-P}+PAx_{I-P}\\
\end{array} \right) \in \hat{A}, 
\]
thus $S=A_{\vert (I-P)\mathcal{K}}$. 

By definition $\mathcal{R}=\left( (I-P)\mathcal{K} \right)^{[\perp]} =P\mathcal{K}$ and $\hat{A}_{\infty}=\mathcal{\tilde{R}}$.

$S$ is a closed symmetric relation in the Pontryagin space $\mathcal{K}$ because it is the intersection of such relations $A$ and $\hat{A}$. $S$ is a bounded operator as a restriction of bounded operator $A$. This proves (iii). 

(iv) Now when we know $S$, we can find $S^{+}$ by definition. It is as claimed in (iv). 

(v) By solving equation $(S^{+}-z)\left( {\begin{array}{*{20}c}
x_{I-P}\\
x_{P}\\
\end{array} } \right)=\left( {\begin{array}{*{20}c}
0\\
0\\
\end{array} } \right)$, i.e., by solving equation
\[ \left( {\begin{array}{*{20}c}
\tilde{A} - z& (I-P)A_{\vert P\mathcal{K}}\\
(I-P)\mathcal{K} \times P\mathcal{K} & P \mathcal{K} \times P\mathcal{K}-z\\
\end{array} } \right)\left( {\begin{array}{*{20}c}
x_{I-P}\\
x_{P}\\
\end{array} } \right)=\left( {\begin{array}{*{20}c}
0\\
0\\
\end{array} } \right)
\]
we obtain 
\[
\mathcal{R}_{z}=\left\lbrace \left( \begin{array}{*{20}c}
-(\tilde{A}-z)^{-1}(I-P)APx_{P}\\
x_{P}\\
\end{array} \right) : x_{P} \in P\mathcal{K} \right\rbrace.
\]
According to \cite[Theorem 4]{B}, the function $\hat{Q}_{2}\left( z \right):=\tilde{\Gamma}_{2}^{+}( \tilde{A}-z )^{-1}\tilde{\Gamma}_{2}\in \mathcal{N}_{\kappa_{2}}\left( \mathcal{H} \right)$, with $\tilde{\Gamma}_{2}:=\left( I-P \right)A\tilde{\Gamma}\left( \tilde{\Gamma}^{+}\tilde{\Gamma} \right)^{-1}$, has $\kappa_{2}$ negative squares, where $\kappa_{2}$ is the negative index of $(I-P)\mathcal{K}$. Then
\begin{equation}
\label{eq330}
(I-P)\mathcal{K}=c.l.s. \left \lbrace ( \tilde{A}-z )^{-1}\tilde{\Gamma}_{2} \mathcal{H}, \, z \in \rho(\tilde{A}) \right \rbrace.
\end{equation}
It is easy to verify
\[
(\tilde{A}-z)^{-1}(I-P)AP\mathcal{K}=( \tilde{A}-z )^{-1}\tilde{\Gamma}_{2}\mathcal{H}=(\tilde{A}-z)^{-1}(I-P)A\tilde{\Gamma }(\tilde{\Gamma }^{+} \tilde{\Gamma })^{-1}\mathcal{H}.
\]
According to (\ref{eq330}) we have
\[
\left( \begin{array}{*{20}c}
f_{I-P}\\
f_{P}\\
\end{array} \right) [\bot]\left( \begin{array}{*{20}c}
-(\tilde{A}-z)^{-1}(I-P)APx_{P}\\
x_{P}\\
\end{array} \right), \forall z \in \rho (A) \Rightarrow  \left( \begin{array}{*{20}c}
f_{I-P}\\
f_{P}\\
\end{array} \right) = \left( \begin{array}{*{20}c}
0\\
0\\
\end{array} \right).
\]
This further means
\[
\mathcal{K}=c.l.s. \left \lbrace \mathcal{R}_{z}: z \in \rho(A) \right \rbrace.
\]
Hence, $S=A_{\vert (I-P)\mathcal{K}}$ is a simple operator in $\mathcal{K}$. 

(vi) If $\tilde{\Gamma}$ is one-to-one, then according to (\ref{eq116}), $\ker \Gamma_{z}=\lbrace 0 \rbrace, \forall z\in \rho(A)$. According to Proposition \ref{proposition22} (i), the function $Q$ is strict. According to Theorem \ref{theorem22} (b), $Q$ is the Weyl function of $A$ corresponding to the boundary triple $\Pi =(\mathcal{H},\Gamma_{0}, \Gamma_{1})$ that satisfies $A=\ker \Gamma_{0}$. The second claim of (vi) follows from Theorem \ref{theorem22} (c). 

The claim $A\dot{+}\mathcal{\hat{R}} = S^{+}$ we can see by comparing elements of the two relations. Indeed, for an arbitrary $f=f_{I-P} +f_{P}\in \mathcal{K}$,
\[
\left\lbrace \left( {\begin{array}{*{20}c}
f_{I-P} +f_{P}\\
\tilde{A}f_{I-P}+(I-P)Af_{P}+PAf_{I-P} +PAf_{P}+P\mathcal{K}\\
\end{array} }\right)\right\rbrace =
\]
\[= \left\lbrace \left( {\begin{array}{*{20}c}
f_{I-P} +f_{P}\\
\tilde{A}f_{I-P}+(I-P)Af_{P}+P\mathcal{K}+P\mathcal{K}\\
\end{array} }\right)\right\rbrace
\]
obviously holds, where we use claim (iv) for $S^{+}$ on the right hand side of the equation.
\end{proof}

Recall that an extension $\tilde{S}\in Ext\, S$ is $\mathcal{R}$-\textit{regular} if $\tilde{S}\hat{+}\hat{\mathcal{R}}$ is a closed linear relation in $\mathcal{K} \times \mathcal{K}$, see \cite [Definition 3.1] {D1}. 

\begin{cor}\label{corollary310} Let $Q\in N_{\kappa }(\mathcal{H})$ be a strict function that satisfies the conditions of Theorem \ref{theorem32}. Then $A$, $\hat A$, and $S^{+}$ are $\mathcal{R}$-regular extensions of $S$.
\end{cor}
\begin{proof} The extension $A$ is $\mathcal{R}$-regular because, according to Theorem \ref{theorem32} (vi), $S^{+}=A\dot{+}\mathcal{\hat{R}}$ and it is a closed relation in $\mathcal{K} \times \mathcal{K}$. 

From $\hat{A}=S \dot{+} \mathcal{\hat{R}}$ and $\mathcal{\hat{R}}\hat{+} \mathcal{\hat{R}}= \mathcal{\hat{R}}$, it follows that $\hat{A}=\hat{A}\hat{+}\mathcal{\hat{R}}$. Since $\hat{A}$ is closed, it is the $\mathcal{R}$-regular extension of $S$. By the same token, $S^{+}$ is $\mathcal{R}$-regular. 
\end{proof}

\section{Examples}\label{s8} 


In the following examples we will show how to use results from sections \ref{s4} and \ref{s6} to find a closed symmetric operator $S$ and a reduction operator $\Gamma$ for a given generalized Nevanlinna function $Q$ so that $Q$ becomes the Weyl function related to $S$ and $\Gamma$. We will also express $S$ and $S^{+}$ in terms of the representing operator $A$ of the function Q. 

\begin{exm}\label{example46} Given function the $Q\left( z \right):=-\frac{1}{z}$ , $Q\in N_{0}(\mathbb{C})$. Find the corresponding symmetric linear realton $S$, $S^{+}$ and the triple $\Pi=(\mathbb{C},\Gamma_{0}, \Gamma_{1})$.
\end{exm}

This function is holomorphic at $\infty $ and 
\[
Q^{'}\left( \infty \right):={\lim \limits_{z \to \infty}{zQ(z)}}=-I_{\mathbb{C}}
\]
is a boundedly invertible operator, i.e. the conditions of Theorem \ref{theorem32} are satisfied. It is also easy to verify that $Q$ is a strict function in $\mathcal{D}(Q)$. According to Lemma \ref{lemma14}, the minimal representation of $Q$ is of the form
\[
Q\left( z \right)=\tilde{\Gamma}^{+}\left( A-z \right)^{-1}\tilde{\Gamma} , z\in 
\rho\left( A \right),
\]
where $A$ is a bounded operator, and $Q^{'}\left( \infty \right)=-\tilde{\Gamma}^{+}\tilde{\Gamma}=-I_{\mathbb{C}}=\left( -1 \right)\in \mathbb{C}^{1\times 1}$. 

We know, and it is easy to verify, that in the representation of the function $Q\left( z \right):=-\frac{1}{z}$, the minimal state space is $\mathcal K=\mathbb{C}$, the representing operator is 
\[
A=(0)=\left\{ \left( 
{\begin{array}{*{20}c}
f\\
0\\
\end{array} } \right):f\in \mathbb{C} \right\}\subseteq \mathbb{C}^{2},
\]
the resolvent is $\left( A-z \right)^{-1}=-\frac{1}{z}I_{\mathbb{C}}$, and $\tilde{\Gamma}^{+}=\tilde{\Gamma} =\left( 1 \right)\in \mathbb{C}^{1\times 1}$ holds. According to (\ref{eq118}), $P=I_{\mathbb{C}}$. Because $P \mathcal{K}=\mathcal{K}$, according to Theorem \ref{theorem32}, $S=A_{\vert (I-P)\mathcal{K}} \cap \hat{A}=\left\lbrace \left( 
{\begin{array}{*{20}c}
0\\
0\\
\end{array} } \right)\right\rbrace $. Then according to Theorem \ref{theorem32} (v), $\mathcal{R}_{z}=P\mathcal{K}=\mathcal{K}$. Because $\tilde{\Gamma}$ is a one-to-one operator, according to Theorem \ref{theorem32} (vi), $Q\left( z \right):=-\frac{1}{z}$ is the Weyl function associated with $S$ and $A$. 

We also know that in the same state space $\mathcal{K}=\mathbb{C}$, there exists a linear relation $\hat{A}$ that minimally represents $\hat{Q}\left( z \right)=-Q^{-1}\left( z \right)=zI_{\mathbb{C}}$, and $\hat{\mathcal{R}}=\left( \left\{ 0 \right\}\times \mathbb{C} \right)\subseteq \mathbb{C}^{2}$. According to Theorem \ref{theorem32} (iii), $\hat{A}=\tilde{A}\dot{[+]}\hat{\mathcal R}=\hat{\mathcal R}$. 

Then, according to Theorem \ref{theorem22} (c) (ii), $ S^{+}=A\hat{+}\hat{A}=\mathbb{C}^{2}$.

Now we need to define the reduction operator $\Gamma =\left( 
{\begin{array}{*{20}c}
\Gamma_{0}\\
\Gamma_{1}\\
\end{array} } \right):\, S^{+} \to \mathcal H^{2}$ that will satisfy identity (\ref{eq124}) and 
\[
A=\ker \Gamma_{0}\, \wedge \hat{A}=\ker \Gamma_{1}.
\]
Because, $M=Q\in N_{0}(\mathbb{C})$, the space $\mathcal{K}=\mathbb{C}$ is endowed with the usual definite scalar product. We can easily verify that the reduction operator that satisfies the above condition is defined by
\[
\Gamma \left( 
{\begin{array}{*{20}c}
f\\
f^{'}\\
\end{array} } \right)=\left( 
{\begin{array}{*{20}c}
f^{'}\\
-f\\
\end{array} } \right).
\] \hfill $\square$ 

\begin{exm}\label{example48} In Example \ref{example22} we derived a strict part $\tilde{Q}(z)=z$ from a non-strict matrix Nevanlinna function. Because the strict part remains a Nevanlinna function and it becomes a strict function, according to Theorem \ref{theorem22} (b) there exist a reduction operator $\Gamma$ and a boundary triple $\Pi$ that correspond to $\tilde{Q}(z)=z$.
\end{exm} 
To accomplish this task, we can use results of Example \ref{example46}, because $-\tilde{Q}(z)^{-1} = -\frac{1}{z}$. This means that $\Gamma_{0}$ and $\Gamma_{1}$ exchange roles, i.e., in this example 
\[
\Gamma \left( {\begin{array}{*{20}c}
f\\
f^{'}\\
\end{array} } \right):=\left( {\begin{array}{*{20}c}
f\\
f^{'}\\
\end{array} } \right).
\] 
Therefore, now we have 
\[
A=\ker\, \Gamma_{0}= \left\lbrace  \left( {\begin{array}{*{20}c}
0\\
f^{'}\\
\end{array} } \right) : f^{'} \in \tilde{\mathcal{H}}\right\rbrace  \wedge 
\hat{A}=\ker\, \Gamma_{1}=\left\lbrace  \left( {\begin{array}{*{20}c}
f\\
0\\
\end{array} } \right) : f \in \tilde{\mathcal{H}} \right\rbrace,
\]
where $ \tilde{\mathcal{H}}=\mathbb{C}$. Then $S=A \cap \hat{A}=\left\lbrace  \left( {\begin{array}{*{20}c}
0\\
0\\
\end{array} } \right) \right\rbrace , S^{+}=\tilde{\mathcal{H}}^{2}$. 
Obviously $ \ker \left( S^{+}-zI\right)=\tilde{\mathcal{H}}$. This implies $\hat{\mathcal{R}}_{z}\left( S^{+}\right)=\left\lbrace  \left( {\begin{array}{*{20}c}
f\\
zf\\
\end{array} } \right) : f \in \tilde{\mathcal{H}} \right\rbrace$. Thus
\[
\Gamma_{0} \left( {\begin{array}{*{20}c}
f\\
zf\\
\end{array} } \right)=f \wedge \Gamma_{1} \left( {\begin{array}{*{20}c}
f\\
zf\\
\end{array} } \right)=zf.
\]
By the definition of the Weyl function, see (\ref{eq126}), it follows that $ \tilde{Q}(z)=z$, i.e. $\tilde{Q}(z)$ is indeed the Weyl function corresponding to the reduction operator $\Gamma$.\hfill $\square$ 

Note that in \cite[Example 2.4.2]{BHS}, the authors start from the symmetric relation $S$ and the reduction operator $\Gamma$ to find the corresponding Weyl function $M$, while in this example we do the converse work, we start from the strict part $\tilde{Q}$ to find $\Gamma$ and $S$. At the end we verified that $\tilde{Q}$ is indeed the Weyl function corresponding to those $\Gamma$ and $S$. 

In the following example, we will show how to use Theorem \ref{theorem32} to find linear relations $S$, $\hat{A}$ and $S^{+}$ for a given function $Q$. 

\begin{exm}\label{example410}  Given the function
\[
Q\left( z \right)=\left( {\begin{array}{*{20}c}
\frac{-(1+z)}{z^{2}} & \frac{1}{z}\\[1mm]
\frac{1}{z} & \frac{1}{1+z}\\[1mm]
\end{array} } \right)\in N_{2}(\mathbb{C}^{2})
\]
and its operator representation
\[
Q\left( z \right)=\tilde{\Gamma}^{+}\left( A\mathrm{-}z \right)^{\mathrm{-1}}\tilde{\Gamma},
\]
where the fundamental symmetry $J$, and operators $A$, $\Gamma$ and $\Gamma^{+}$ are, respectively:
\[
J=\left( {\begin{array}{*{20}c}
0 & 1 & 0\\
1 & 0 & 0\\
0 & 0 & -1\\
\end{array} } \right), A=\left( {\begin{array}{*{20}c}
0 & 1 & 0\\
0 & 0 & 0\\
0 & 0 & -1\\
\end{array} } \right), \tilde{\Gamma} =\left( {\begin{array}{*{20}c}
0.5 & -1\\
1 & 0\\
0 & -1\\
\end{array} } \right), \tilde{\Gamma}^{+}=\left( 
{\begin{array}{*{20}c}
1 & 0.5 & 0\\
0 & -1 & 1\\
\end{array} } \right),
\]
our task is to find linear relations $S$, $\hat{A}$ and $S^{+}$. 
\end{exm}
It is easy to verify that the function $Q$ satisfies the conditions of Theorem \ref{theorem32}. Indeed, the limit (\ref{eq114}) gives
\[
\tilde{\Gamma}^{+}\tilde{\Gamma} =\left( {\begin{array}{*{20}c}
1 & -1\\
-1 & -1\\
\end{array} } \right), \left( \tilde{\Gamma}^{+}\tilde{\Gamma} \right)^{-1}=\left( 
{\begin{array}{*{20}c}
0.5 & -0.5\\
-0.5 & -0.5\\
\end{array} } \right).
\]
Then, by means of formula (\ref{eq118}), we get 
\[
P=\left( {\begin{array}{*{20}c}
0.75 & 0.125 & 0.25\\
0.5 & 0.75 & -0.5\\
0.5 & -0.25 & 0.5\\
\end{array} } \right), I-P=\, \left( {\begin{array}{*{20}c}
0.25 & -0.125 & -0.25\\
-0.5 & 0.25 & 0.5\\
-0.5 & 0.25 & 0.5\\
\end{array} } \right).
\]
According to Theorem \ref{theorem32} (iii), we can find $S$:
\[
S= A\left( I-P\right)=\left( {\begin{array}{*{20}c}
-0.5 & 0.25 & 0.5\\
0 & 0 & 0\\
0.5 & -0.25 & -0.5\\
\end{array} } \right).
\] 
\[
\tilde{A}:=\left( I-P\right) A\left( I-P\right)=\left( {\begin{array}{*{20}c}
-0.25 & 0.125 & 0.25\\
0.5 & -0.25 & -0.5\\
0.5 & -0.25 & -0.5\\
\end{array} } \right)=-(I-P).
\] 
By solving equation $Px=x$ and then using the fact $(I-P)\mathcal{K}[\perp]P\mathcal{K}$, we obtain
\[
(I-P) \mathcal K =l.s.\left\lbrace \left( {\begin{array}{*{20}c}
-1\\
2\\
2\\
\end{array} } \right)\right\rbrace; \,\, P \mathcal{K} =l.s.\left\lbrace \left( {\begin{array}{*{20}c}
3\\
2\\
2\\
\end{array} } \right),
 \left( {\begin{array}{*{20}c}
1 \\
0\\
1\\
\end{array} } \right)
\right\rbrace.
\]
According to Theorem \ref{theorem32} (ii) we have 
\[
\hat{A}=\tilde{A} \dot{[+]}\hat{\mathcal{R}}=-I_{I-P}\dot{[+]}\left( \lbrace 0 \rbrace \times P \mathcal K\right).
\] 
The equivalent, developed form of the linear relation $\hat{A}$ is:  
\[
\hat{A}\left( {\begin{array}{*{20}c}
f_{1}\\
f_{2}\\
f_{3}\\
\end{array} } \right)=\left(\frac{f_{1}}{4}-\frac{f_{2}}{8}-\frac{f_{3}}{4} \right)
\left( {\begin{array}{*{20}c}
-1\\
2\\
2\\
\end{array} } \right)+
c_{1}\left( {\begin{array}{*{20}c}
3\\
2\\
2\\
\end{array} } \right)+
c_{2}\left( {\begin{array}{*{20}c}
1 \\
0\\
1\\
\end{array} } \right),
\]
where $f=\left( {\begin{array}{*{20}c}
f_{1}\\
f_{2}\\
f_{3}\\
\end{array} } \right) \in \mathcal{K}=\mathbb{C}^{3}$, and $ c_{i} \in \mathbb{C},i=1, 2 $, are arbitrary constants.

The easiest way to obtain the developed form of $S^{+}$ is to use Theorem \ref{theorem32} (vi) representation $S^{+}=A \dot{+} \mathcal{\hat{R}}$. We get
\[
S^{+}f=\left( {\begin{array}{*{20}c}
0 & 1 & 0\\
0 & 0 & 0\\
0 & 0 & -1\\
\end{array} } \right)\left( {\begin{array}{*{20}c}
f_{1}\\
f_{2}\\
f_{3}\\
\end{array} } \right) + P \mathcal{K}=\left( {\begin{array}{*{20}c}
f_{2}\\
0\\
-f_{3}\\
\end{array} } \right) + c_{1}\left( {\begin{array}{*{20}c}
3\\
2\\
2\\
\end{array} } \right)+
c_{2}\left( {\begin{array}{*{20}c}
1 \\
0\\
1\\
\end{array} } \right),
\]
where $f$ and $c_{i} \in \mathbb{C}, \,i=1,2$, are as before. \hfill $\square$
\\







\vspace{.5cm}
\begin{footnotesize}
	\noindent
	\begin{tabular}{l}
\\ \\ \\ \\ \\ \\ \\ \\ \\ \\ \\
\end{tabular}
	\hfill
	\begin{tabular}{l}
Muhamed Borogovac\\
Boston Mutual Life\\
Actuarial Department\\
120 Royall St. Canton, MA 02021, USA\\
e-mail: {muhamed.borogovac@gmail.com }\\
\end{tabular}
\end{footnotesize}
\label{LastPage}
\end{document}